\documentclass[12pt]{article}
\usepackage{latexsym,color,amsmath,amsthm,amssymb,amscd,amsfonts}

\newtheorem{lemma}{Lemma}[section]
\newtheorem{proposition}[lemma]{Proposition}
\newtheorem{remark}[lemma]{Remark}
\newtheorem{theorem}[lemma]{Theorem}
\newtheorem{definition}[lemma]{Definition}

\newtheorem{claim}[lemma]{Claim}

\newtheorem*{remark*}{Remark}

\newcommand{\Id}{{{\mathchoice {\rm 1\mskip-4mu l} {\rm 1\mskip-4mu l}
      {\rm 1\mskip-4.5mu l} {\rm 1\mskip-5mu l}}}}

\begin{document}
\title {On the extremality of Hofer's metric on the group of Hamiltonian
diffeomorphisms}
\author{Yaron Ostrover and Roy Wagner }

\maketitle
\begin{abstract}
\noindent Let $M$ be a closed symplectic manifold, and let
$\|\cdot\|$ be a norm on the space of all smooth functions on $M$,
which are zero-mean normalized with respect to the canonical
volume form. We show that if $\|\cdot\| \leq C\|\cdot\|_{\infty}$,
and $\|\cdot\|$ is invariant under the action of Hamiltonian
diffeomorphisms, then it is also invariant under all volume
preserving diffeomorphisms. We also prove that if $\|\cdot\|$ is,
additionally, not equivalent to $\|\cdot\|_{\infty}$, then the
induced Finsler metric on the group Ham$(M,\omega)$ of Hamiltonian
diffeomorphisms on $M$ vanishes identically. These results provide
partial answers to questions raised by Eliashberg and Polterovich
in~\cite{EP}. Both results rely on an extension of $\|\cdot\|$ to
the space of essentially bounded measurable functions, which is
invariant under all measure preserving bijections.
\end{abstract}

\section{Introduction and Results}

Let $(M,\omega)$ be a closed connected symplectic manifold of
dimension $2n$. Denote by ${\cal A}$ the space of all smooth
functions on $M$ which are zero-mean normalized with respect to
the canonical volume form $\omega^n$. The main object of our study
is the infinite-dimensional Lie group Ham$(M,\omega)$ of
Hamiltonian diffeomorphisms of $M$. We refer the readers
to~\cite{HZ},~\cite{MS} and~\cite{P1} for symplectic preliminaries
and further discussions on the group of Hamiltonian
diffeomorphisms.

\noindent It is well known that the Lie algebra of
Ham$(M,\omega)$, that is the space of all Hamiltonian vector
fields, can be identified with the space ${\cal A}$. Moreover, the
adjoint action of Ham$(M,\omega)$ on its Lie algebra ${\cal A}$ is
the standard action of diffeomorphisms on functions. The choice of
any norm $\| \cdot \|$ on ${\cal A}$ gives rise to a
pseudo-distance function on Ham$(M,\omega)$ in the following way:
we define the length of a path $\alpha : [0,1] \rightarrow {\rm
Ham}(M,\omega)$ as
$$ {\rm length}\{ {\alpha}\} = \int_0^1 \| \dot \alpha_t \| dt =
\int_0^1 \| F_t \| dt ,$$ where $F_t(x)=F(t,x)$ is the Hamiltonian
function generating the path $\alpha$. This is the usual
definition of Finsler length. The distance between two Hamiltonian
diffeomorphisms is given by
$$ \rho(\psi,\varphi) = \inf {\rm length} { \{ \alpha \} },$$ where the infimum
is taken over all Hamiltonian paths $\alpha$ connecting $\psi$
and $\varphi$. It is not hard to check that $\rho$ is
non-negative, symmetric and satisfies the triangle inequality.
Moreover, a norm on ${\cal A}$ which is invariant under the
adjoint action yields a bi-invariant pseudo-distance function,
i.e.
 $\rho(\psi,\varphi) = \rho(\theta \, \psi,\theta \, \varphi) = \rho(\psi \, \theta ,\varphi \, \theta)$ for every
 $\psi,  \varphi,  \theta \in {\rm Ham} (M,\omega)$. From now on we will deal only with such norms and we will refer to
$\rho$ as the pseudo-distance generated by the norm $\| \cdot \|$.

\noindent It is highly non-trivial to check whether such a norm is
non-degenerate, that is $\rho({\mathbb \Id},\psi) > 0$ for $\psi
\neq {\mathbb \Id}$. In fact, for compact symplectic manifolds, a
bi-invariant pseudo-metric $\rho$ on Ham$(M,\omega)$ is either a
genuine metric or identically zero. This is an immediate corollary
of a well known theorem by Banyaga~\cite{B}, which states that
Ham$(M,\omega)$ is a simple group, combined with the fact that the
null-set $${\rm null}(\rho) = \{ \psi \in {\rm Ham}(M,\omega) \ |
\ \rho(\Id,\psi) = 0 \}$$ is a normal subgroup of Ham$(M,\omega)$.

\noindent A distinguished result by Hofer~\cite{H} states that the
$L_{\infty}$ norm $\| \cdot \|_{\infty}$ on ${\cal A}$ gives rise
to a genuine distance function on Ham$(M,\omega)$. This was
discovered and proved by Hofer for the case of ${\mathbb R}^{2n}$,
then generalized by Polterovich~\cite{P2} to some larger class of
symplectic manifolds, and finally proven in full generality by
Lalonde and McDuff in~\cite{LM}. The above mentioned distance
function is known as Hofer's metric and has been intensively
studied since its discovery (see
e.g.~\cite{HZ},~\cite{MS},~\cite{P1}).
 We also refer the reader to Oh's
paper~\cite{O} for another approach to the non-degeneracy of
Hofer's metric. In the opposite direction, Eliashberg and
Polterovich showed in~\cite{EP} that for $1 \leq p < \infty$, the
pseudo-distances on Ham$(M,\omega)$ which correspond to the $L_p$
norms on ${\cal A}$ vanishes identically. Thus, the following
question arises from~\cite{EP} and~\cite{P1}:

\noindent {\bf Question:} What are the invariant norms on ${\cal
A}$, and which of them give rise to genuine bi-invariant metrics
on Ham$(M,\omega)$?

\noindent Our main contributions towards answering this question
are

\begin{theorem} \label{Ham-invariant-implies-measure-invarinat}
Any Ham$(M,\omega)$-invariant norm $\| \cdot \|$ on ${\cal A}$
which is dominated from above by $ \| \cdot \|_{\infty}$ is
invariant under all measure preserving diffeomorphisms on $M$.
\end{theorem}

\begin{theorem} \label{The_main_theorem} Let $ \| \cdot \|$
be a Ham$(M,\omega)$-invariant norm on ${\cal A}$ such that $\|
\cdot \| \leq C\| \cdot \|_{\infty}$ for some constant $C$, but
the two norms are not equivalent. Then the associated
pseudo-distance function $\rho$ on Ham$(M,\omega)$ vanishes
identically.
\end{theorem}

\noindent Here, two norms are said to be \emph{equivalent}, if
each dominates the other up to a multiplicative constant.

\noindent The next result is a strengthened formulation of
Theorem~\ref{Ham-invariant-implies-measure-invarinat} and a key
ingredient in the proof of Theorem~\ref{The_main_theorem}. As the
discussion below explains, it also bears on the question of
classifying Ham$(M,\omega)$-invariant norms.

\begin{theorem}\label{extension-to-L-infinity}
Let $ \| \cdot \|$ be a Ham$(M,\omega)$-invariant norm on ${\cal
A}$ such that $\| \cdot \| \leq C \| \cdot \|_{\infty}$ for some
constant $C$. Then $\| \cdot \|$ can be extended to a (semi) norm
$||| \cdot ||| \leq C\| \cdot \|_{\infty}$ on $L_{\infty}(M)$,
which is invariant under all measure preserving bijections on $M$.
\end{theorem}

\noindent The formulation of the theorem states only what is
necessary for the proofs of
Theorems~\ref{Ham-invariant-implies-measure-invarinat}
and~\ref{The_main_theorem}. In fact we know more about
$|||\cdot|||$. First, $|||\cdot|||$ is a norm, rather than just a
seminorm (namely, $|||\cdot|||$ does not vanish on non-zero
functions). Second, $|||\cdot|||$ coincides (where it is defined)
with the completion of $||\cdot||$, which in turn can be viewed as
a dense subspace of $L_1(M)$ equipped with a norm invariant under
measure preserving bijections. The argument for the first claim is
briefly sketched in Remark~\ref{non-degenracy}, and that for the
second claim is outlined in the final section. The final section
also refers to literature concerning the classification of such
norms, and indicates their possible pathologies.

\noindent {\bf Structure of the paper:}

\noindent The next section contains a fairly detailed outline of
the proofs of our main theorems, stressing the main ingredients
involved. The following two sections present complete proofs of
Theorem~\ref{extension-to-L-infinity} and
Theorem~\ref{The_main_theorem} respectively. Section $5$ contains
proofs of some lemmas. The last section contains a sketchy
treatment of additional properties of the norm $|||\cdot|||$,
together with some references concerning the classification of
such norms.

\section{Outline of the Proofs}

\noindent As explained in the introduction, the degeneracy of the
pseudo-distance function $\rho$ (Theorem~\ref{The_main_theorem}),
is proved in~\cite{EP} for $L_p$ norms, $1 \leq p < \infty$. The
only property of $L_p$ actually used in that proof is, roughly
speaking, that uniformly bounded functions with small support have
small norm. More precisely, in section 4 we reproduce an argument
from from~\cite{EP} to show that the proof of
Theorem~\ref{The_main_theorem} can be reduced to the following
claim

\begin{claim} \label{claim-smallsupp->smallnorm} If $\sup\{\|F_n\|_{\infty}\} < \infty$ and ${\rm
Vol} \bigl ({\rm Support}(F_n) \bigr) \rightarrow 0$, then $\| F_n
\| \rightarrow 0$.
\end{claim}

\noindent Therefore, our main task is to prove this property for
any norm which satisfies the requirements of
Theorem~\ref{The_main_theorem}. As will be explained below,
Theorem~\ref{extension-to-L-infinity} allows us to carry out the
proof of this claim in a more amenable setting.

\noindent A natural approach to
Claim~\ref{claim-smallsupp->smallnorm} would be to consider
characteristic functions with small-measure support first, then
make the standard move to step functions, and conclude with any
smooth bounded function with small-measure support. The obvious
obstacle is that characteristic functions are not smooth, and are
therefore outside our space. Here one may choose to approximate
them by smooth functions and work from there. We chose, however,
to extend our setting so as to include genuine characteristic
functions. This is where Theorem~\ref{extension-to-L-infinity}
comes in. We will interrupt the discussion on the proof of
Claim~\ref{claim-smallsupp->smallnorm} to discuss the proof of
Theorem~\ref{extension-to-L-infinity}.

\noindent Recall that our aim in
Theorem~\ref{extension-to-L-infinity} is to extend the given norm
$\| \cdot \|$ to $L_{\infty}(M)$. For this purpose, we first
extend our norm to all smooth functions, with average not
necessarily zero (since this adds just one dimension to our
original space of functions, any two extensions are equivalent).
Next, we take advantage of the fact that $C^{\infty}(M)$ is dense
in $L_{\infty}(M)$ with respect to the topology of convergence in
measure. We define
$$ |||F||| = \inf \{ \liminf_{n \rightarrow \infty} ||F_n|| \}, $$ where the infimum is
taken over all sequences $\{F_n\}$ of uniformly bounded smooth
functions which converge in measure to $F$.

\noindent Such constructions occur occasionally in functional
analysis, for instance in the extension of the Riemann integral
from continuous to semi-continuous functions (using pointwise
convergence from above/below), and in the extension of operator
norm from finite-rank operators on a Banach space to approximable
operators (using uniform convergence on compacts). However, we
are not aware of any similar construction which relies on
convergence in measure.

\noindent We study $|||\cdot|||$ in section $3$. First we confirm
that $||| \cdot |||$ is a semi-norm on $L_{\infty}(M)$ which is
dominated from above by $\| \cdot \|_{\infty}$. We then go on to
prove the non-trivial properties of $||| \cdot |||$: it coincides
with $\| \cdot \|$ on smooth functions, and is invariant under
measure preserving bijections. Formally:

\begin{claim} \label{claim-coincidness}
For every $F \in {\cal A}$ we have $\|F\|=|||F|||$.
\end{claim}

\begin{claim} \label{claim-measure-preserving-invariant}
For every $F \in L_{\infty}(M)$ and every measure preserving
bijection $\varphi$ on $M$ we have $$||| F \circ \varphi ||| =
|||F|||$$
\end{claim}

\noindent In order to prove this second property, recall that our
original norm $\| \cdot \|$ is already invariant under Hamiltonian
diffeomorphisms. To extend the invariance we invoke Katok's
``Basic Lemma" from~\cite{K}, which effectively allows to
approximate in measure any measure preserving bijection by a
Hamiltonian diffeomorphism. More precisely, fix an arbitrary
Riemannian metric $d$ on $M$. We claim
\begin{lemma} \label{approximation in measure lemma}
For every measure preserving bijection $\varphi$ of $M$ (not
necessarily continuous) and every $\varepsilon > 0$, there exists
a Hamiltonian diffeomorphism $g$ on $M$ which
$\varepsilon$-approximates $\varphi$ in measure, namely
$$ {\rm Vol} \bigl (\{x \in M \, ; \, d(\varphi(x),g(x)) > \varepsilon\} \bigr ) < \varepsilon $$
\end{lemma}

\noindent This result is of course independent of the specific
Riemannian structure chosen. We postpone the proof of the lemma
to the last section of this paper. The proof of
Claim~\ref{claim-measure-preserving-invariant} follows easily from
Lemma~\ref{approximation in measure lemma} and the definition of
$||| \cdot |||$. Claim~\ref{claim-coincidness} and
Claim~\ref{claim-measure-preserving-invariant} conclude the proof
of Theorem~\ref{extension-to-L-infinity}.

\noindent With a measure-preserving-bijection-invariant extension
of $||\cdot||$ at our disposal, let's return to the proof of
Claim~\ref{claim-smallsupp->smallnorm}. Note that
Claim~\ref{claim-coincidness} implies that it is sufficient to
prove Claim~\ref{claim-smallsupp->smallnorm} for the norm $|||
\cdot |||$. The rest of this section is devoted to this issue.

\noindent Our argument depends on the fact, inspired by an
argument from~\cite{P3}, that an operator, which performs
piecewise averaging on functions, is bounded. More precisely,
relying on the fact that $||| \cdot |||$ is invariant under
measure preserving bijections, we prove that

\begin{lemma}[Piecewise-Averaging property] \label{claim-picewise-averaging}   For every continuous $F$ and every measurable partition
$\{S_i\}$ of $M$, we have $$|||\sum_i \langle F \rangle_{S_i}
{\mathbb \Id}_{S_i}||| \leq |||F|||,$$ where $\langle F
\rangle_{S_i} = {\frac 1 { {\rm Vol}(S_i)}} \int_{S_i}F \omega^n$
 denotes the average of $F$ over $S_i$.
\end{lemma}

\noindent The proof of the lemma is postponed to the last section.
Let us now explain how this property serves to prove
Claim~\ref{claim-smallsupp->smallnorm}. Fix $\varepsilon
> 0$. The hypothesis of Theorem~\ref{The_main_theorem} provides us
with smooth functions $F$ such that $\|F\|_{\infty}=1$ while
$\|F\| = |||F ||| \leq \varepsilon$. Partition $M$ into $A$ and
$A^c = M \setminus A$, where $A$ is a small enough neighborhood of
the maximum of $F$, such that $ || {\mathbb \Id}_A -
 \langle F \rangle_A {\mathbb \Id}_A + \langle F \rangle_{A^c} {\mathbb
\Id}_{A^c}||_{\infty} < \varepsilon.$ Next, it follows from
Lemma~\ref{claim-picewise-averaging}, the fact that $||| \cdot
|||$ is dominated from above by $\| \cdot \|_{\infty}$, and the
triangle inequality that
$$ |||{\mathbb \Id}_A ||| \leq || {\mathbb \Id}_A  -
 \langle F \rangle_A {\mathbb \Id}_A + \langle F \rangle_{A^c} {\mathbb
\Id}_{A^c}||_{\infty} + ||| \langle F \rangle_A {\mathbb \Id}_A +
\langle F \rangle_{A^c} {\mathbb \Id}_{A^c}||| \leq \varepsilon +
||| F||| \leq 2 \varepsilon$$

\noindent Since $|||\cdot|||$ is invariant under measure
preserving bijections, this applies to every set $B$ with the same
measure as $A$. Thus, we establish
Claim~\ref{claim-smallsupp->smallnorm} for sequences of
characteristic functions on sets with measure tending to zero. It
is now a simple approximation argument which establishes the Claim
as stated for smooth functions. The complete details are given in
section $4$.

\section{Proof of Theorem~\ref{extension-to-L-infinity}}

\noindent In this section we construct the norm $|||\cdot|||$, and
prove its properties as stated in
Theorem~\ref{extension-to-L-infinity}. The first step towards the
construction of $|||\cdot|||$ is an extension of the given norm to
$C^{\infty}(M)$. Let $C$ such that $\|\cdot\| \leq
C\|\cdot\|_{\infty}$. Endow the space $C^{\infty}(M)$ of all
smooth function on $M$ with the norm $\| \cdot \|'$ defined by
$$ \| F \|' = \inf \bigl \{ \|F_1\| + C\|F_2\|_{\infty} \ ;
 \  F=F_1 + F_2, \ \ F_1 \in {\cal A}, \ \ F_2 \in C^{\infty}(M)
 \bigr \}.$$

\noindent The above definition is just the analytic presentation
of the norm corresponding to the convex hull of the unit ball of
$({\cal A},\|\cdot\|)$ with the unit ball of
$(C^{\infty}(M),\|\cdot\|_{\infty})$, the latter homothetically
shrunk so as to fit inside the former when restricted to ${\cal
A}$. The homogeneity of the new norm is clear. To see that the new
norm satisfies the triangle inequality, let $F=F_1+F_2$ and
$G=G_1+G_2$ such that $\|F_1\| + C\|F_2\|_{\infty} \leq \|F\|' +
\varepsilon$, and $\|G_1\| + C\|G_2\|_{\infty} \leq \|G\|' +
\varepsilon$. Then
\begin{eqnarray*}
\|F+G\|' & \leq & \|F_1+G_1\| + C\|F_2+G_2\|_{\infty}   \\
 & \leq &  (\|F_1\| + C\|F_2\|_{\infty}) +   (\|G_1\| + C\|G_2\|_{\infty})
  \\ & \leq  & \|F\|' + \|G\|' +2\varepsilon.
\end{eqnarray*}

\noindent The new norm is obviously Ham$(M,\omega)$-invariant. To
see that $\|F \|' \leq C\|F\|_{\infty}$, just substitute $F_1=0$
and $F_2=F$ in the definition. To see that $\|\cdot\|' =
\|\cdot\|$ on ${\cal A}$, let $F=F_1+F_2$ where $F_1 \in {\cal A}$
and $F_2 \in C^{\infty}(M)$. Choosing $F_1=F$ and $F_2=0$ proves
that $\|F\|' \leq \|F\|$. For the opposite direction note that
since $F,F_1 \in {\cal A}$, and since $F_2=F-F_1$, the function
$F_2$ must also be in ${\cal A}$. Therefore
$$\|F\|' = \inf_{_{F=F_1+F_2}} \bigl \{ \|F_1\| + C \|F_2\|_{\infty} \bigr \} \geq
\inf_{_{F=F_1+F_2}} \bigl \{ \|F_1\|+\|F_2\| \bigr \} \geq
\inf_{_{F=F_1+F_2}}  \|F_1+F_2\| = \|F\| $$

\noindent Now we are ready to extend our norm to the entire
$L_{\infty}(M)$. Using the same convex-hull trick won't do (it
will fail invariance under measure preserving bijections).
Instead, we will take advantage of the classical fact that any
measurable function can be approximated in measure arbitrarily
well by smooth functions (see e.g.~\cite{R}). We will define a new
functional by taking the least $\|\cdot\|'$ norm among all such
approximations. Formally, we will endow the space $L_{\infty}(M)$
with
$$ ||| F |||  = \inf \left \{ \liminf_{n \rightarrow \infty}
\|F_n\|' \right \},$$ where the infimum is taken over all
sequences of uniformly bounded smooth functions $\{ F_n \} $ which
converge in measure to $F$.

\noindent It is clear that the new functional is homogeneous. To
see that it obeys the triangle inequality, take $\{F_n\}$ and
$\{G_n\}$ which satisfy $\liminf \|F_n\|' \leq |||F||| +
\varepsilon$ and $\liminf \|G_n\|' \leq |||G||| + \varepsilon$.
Then
$$|||F+G||| \leq \liminf \|F_n+G_n\|' \leq \liminf (\|F_n\|' +
\|G_n\|') \leq |||F||| + |||G||| + 2\varepsilon .$$

\noindent To see that the new functional is still bounded by
$C\|\cdot\|_{\infty}$, note that any essentially bounded function
$F$ can be approximated in measure by smooth $F_n$'s with at most
the same essential supremum. Indeed, take any approximation in
measure $F_n$ of $F$, and replace it with ${\rm sign}(F_n) \cdot
(f_n \circ | F_n | )$, where $f_n$ is a good enough smooth
approximation from below of the function $f(s) : {\mathbb R}^+
\rightarrow {\mathbb R}^+$ defined by $f(s) =
\min\{s,\|F\|_{\infty}\}$. Taking such $F_n$'s we get
$$|||F||| \leq
\liminf \|F_n\|' \leq C \liminf \|F_n\|_{\infty} \leq
C\|F\|_{\infty}.$$

\noindent In order to complete the proof of
Theorem~\ref{extension-to-L-infinity} we need the following two
claims.

\begin{claim}\label{extension}
For every $F \in {\cal A}$ we have $\|F\| = \|F\|' = |||F|||$
\end{claim}

\begin{claim}\label{invariance}
For every $F \in L_{\infty}(M)$ and every measure preserving
bijection $\varphi$ on $M$ we have $$||| F \circ \varphi ||| =
|||F|||$$
\end{claim}

\noindent In order to prove the first claim a certain technical
lemma is needed. To state the lemma, fix from now on an arbitrary
Riemannian structure on $M$, and denote by $d$ the corresponding
distance function. Our results are, of course, independent of the
specific Riemannian structure chosen.

\begin{lemma}[Covering Evenly by Many Packings]\label{good_partition_lemma}
 For every $\delta > 0$ and  $\varepsilon > 0$ there exists
a covering  of $M$ by connected open subsets $\{U_i^j\}$, where
$j=1,\ldots,J$ and $i=1,\ldots,L_j$, such that
\renewcommand{\labelenumi}{(\roman{enumi})}
\begin{enumerate}
\item for every fixed $j$, each pair of sets $\{ U_i^{j} \}$ have a positive distance from each other.
\item the diameter of $U_i^j$ with respect to $d$ is at most $\delta$
for all $i$ and $j$.
\item for every $x \in M$, the number of
$j$'s for which $x \notin \cup_i U_i^j$ is at most $\varepsilon
J$.
\end{enumerate}
\renewcommand{\labelenumi}{\arabic{enumi}.}
\end{lemma}

\noindent The proof of the lemma is postponed to the last section
of this paper.

\begin{proof}[{\bf Proof of Claim~\ref{extension}:}]
The restricted equality $\|\cdot\|=\|\cdot\|'$ has been proved
along with the definition of $\|\cdot\|'$ above. Let's prove the
restricted equality $|||\cdot|||=\|\cdot\|'$. By choosing $F_n=F$
for all $n$ in the definition of $|||\cdot|||$, we get $||| \cdot
||| \leq \| \cdot \|'$. In order to show that $\| \cdot \|' \leq
||| \cdot |||$, let $F \in {\cal A}$ and let $\{F_n\}$ be a
sequence of uniformly bounded smooth functions, which converges in
measure to $F$. We need to show that
$$ \liminf_{n \rightarrow \infty} \| F_n \|' \geq  \| F
\|'. $$ For this purpose we will construct a sequence $\{
{\widetilde F_n} \}$ which converges uniformly to $F$, such that $
 \| F_n \|'  \geq \| {\widetilde F_n} \|'$.
Since $\| \cdot \|' \leq C\| \cdot \|_{\infty}$, uniform
convergence implies convergence in $\|\cdot\|'$, and we can
conclude
$$
 \liminf_{n \rightarrow \infty} \| F_n \|' \geq
 \liminf_{n \rightarrow \infty} \| {\widetilde F_n} \|'
 =  \| F \|' .$$

\noindent Let us construct the sequence $\{ {\widetilde F_n} \}$.
Fix $\varepsilon>0$, and let $\delta > 0$ such that every open
neighborhood of diameter $2\delta$ in $M$ can be viewed as a
neighborhood in ${\mathbb R}^{2n}$ such that the original $d$ and
the Euclidian distance are equivalent up to a factor of $2$. Take
a covering $\{U_i^j\}$ of $M$ as in
Lemma~\ref{good_partition_lemma} with the given $\varepsilon$ and
$\delta$. Take $\eta<\delta/6$ such that the $3\eta$-extensions of
any two sets $U_i^j$ with the same $j$ still have a positive
distance between them, and such that
\begin{equation}\label{def-eta}
d(x,y) \leq 2\eta \ \rightarrow \ |F(x) - F(y)|\leq \varepsilon .
\end{equation}
\noindent Set $V_i^j$ to be the $3\eta$-extension of $U_i^j$ with
respect to the distance $d$ on $M$. $\eta$ was chosen such that
$V_i^j$ has diameter at most $2\delta$, and can therefore be
viewed as a neighborhood in ${\mathbb R}^{2n}$ where $d$ and the
Euclidean distance are equivalent up to a factor of $2$. In
particular, any closed Euclidean ball of radius $\eta$ centered
inside $U_i^j$ is contained in $V_i^j$. Denote by $B_{\eta}(x)$
the Euclidean ball of radius $\eta$ around $x$.
Requirement~(\ref{def-eta}) guarantees that
\begin{equation}\label{eta} \left | \langle F \rangle_{B_{\eta}(x)} -  F(x)
\right | \leq \varepsilon .
\end{equation}
\noindent Next, set $n$ such that
$$ {\rm Vol} ( \{x : |F_n(x) - F(x)| > \varepsilon \} ) \ < \
\frac{\varepsilon\cdot
|B_{\eta}|}{\max\{||F_n||_{\infty},||F||_{\infty}\}} , $$ where
$|B_{\eta}|$ is the measure of a Euclidean ball of radius $\eta$.
This is possible since $\{F_n\}$ converges to $F$ in measure, and
since the $F_n$'s are uniformly bounded. This choice of $n$
implies that
\begin{equation}\label{n}
\left | \langle F_n \rangle_{B_{\eta}(x)} - \langle F
\rangle_{B_{\eta}(x)} \right  | \leq 3\varepsilon .
\end{equation}

\noindent By the definition of the integral, and the uniform
continuity of $F_n$, there exist points $\{x^k\}_{k=1}^K \subseteq
B_{\eta}(0)$ such that for every $x \in U_i^j$
$$ \left \lvert {\frac 1 K} \sum_{k=1}^K F_n \left ( x+x^k \right ) -
\langle F_n \rangle_{B_{\eta}(x)} \right \rvert \leq \varepsilon
.$$

\noindent Note that we have established that $V_i^j$ contains the
closure of the $\eta$-extension of $U_i^j$. Thus, using a standard
cut-off argument, we consider Hamiltonian diffeomorphisms
 $g_{i,j}^1,\ldots,g_{i,j}^K$, all
supported inside $V_i^j$, defined by $g_{i,j}^k(x) = x+x^k$ inside
$U_i^j$ and $g_{i,j}^k(x) = x$ outside a small neighborhood of
$U_i^j$. We therefore get for all $x \in U_i^j$
\begin{equation}\label{F_n} \left \lvert {\frac 1 K} \sum_{k=1}^K F_n \left ( g_{i,j}^k(x) \right ) -
\langle F_n \rangle_{B_{\eta}(x)} \right \rvert \leq \varepsilon .
\end{equation}

\noindent Note that for fixed $j$ and $k$, the Hamiltonian
diffeomorphisms $\{g_{i,j}^k\}$ have disjoint supports, and can
therefore be bundled together to form a single diffeomorphism. We
set
$$ {\widetilde {F_n}}(x)  = {\frac 1 {J}} \sum_{j=1}^J \Bigl
({\frac 1 {K}} \sum_{k=1}^K F_n \bigl (\prod_i g_{i,j}^k(x) \bigr
) \Bigr ).$$

\noindent From the triangle inequality and the fact that the norm
$\| \cdot \|'$ is invariant under Hamiltonian diffeomorphisms we
conclude that $\| {\widetilde {F_n}} \|' \leq \|F_n\|'$. Hence, we
need only show that $ \|{\widetilde {F_n}} - F\|_{\infty}
\rightarrow 0$ as $\varepsilon \rightarrow 0$. Indeed,

$$ {\widetilde {F_n}}(x)  =  {\frac 1 {J}}  \biggl ( \sum_{j \in
{\cal J}(x) } \Bigl ( {\frac 1 {K}}\sum_k F_n \bigl (\prod_i
g_{i,j}^k(x) \bigr ) \Bigr ) + \sum_{j \in {\cal J}^c(x) } \Bigl (
{\frac 1 {K}}\sum_k F_n \bigl (\prod_i g_{i,j}^k(x) \bigr ) \Bigr
) \biggr ),
$$

\noindent where ${\cal J}(x) =  \{ j \, | \, x \in \cup_i U_i^j\}
$, ${\cal J}^c(x) =  \{ j \, | \, x \notin \cup_i U_i^j\} $.
 Recall that the third item of
Lemma~\ref{good_partition_lemma} limited the cardinality of ${\cal
J}^c(x)$ to at most $\varepsilon J$ for all $x$. Together with~(\ref{F_n}) this implies that
$$ \Bigl | {\widetilde {F_n}}(x) - {\frac 1 {J}}  \sum_{j=1}^J \bigl (
\langle F_n \rangle_{B_{\eta}(x)} \bigr ) \Bigr | \leq \varepsilon
{\frac { |{\cal J}(x)| } J } + {\frac { |{\cal J}^c(x)| } J }
\cdot 2\max \|F_n\|_{\infty}  \leq \varepsilon + 2\varepsilon
\cdot \max \|F_n\|_{\infty} $$ Together with (\ref{eta}) and
(\ref{n}) we conclude that
\begin{eqnarray*}  \bigl | {\widetilde {F_n}}(x) - F(x) \bigr | & \leq & \Bigl |
{\widetilde {F_n}}(x) - {\frac 1 {J}}  \sum_{j=1}^J \bigl (
\langle F_n \rangle_{B_{\eta}(x)} \bigr ) \Bigr | + \Bigl | {\frac
1 {J}}  \sum_{j=1}^J \bigl ( \langle F_n \rangle_{B_{\eta}(x)}
\bigr ) - {\frac 1 {J}}  \sum_{j=1}^J \bigl ( \langle F
\rangle_{B_{\eta}(x)} \bigr )
 \Bigr | \\ & + & \Bigl | {\frac 1 {J}}  \sum_{j=1}^J \bigl ( \langle
F \rangle_{B_{\eta}(x)} \bigr ) - {\frac 1 {J}}  \sum_{j=1}^J
\bigl ( F(x) \bigr )
 \Bigr | + \Bigl | {\frac 1 {J}}  \sum_{j=1}^J \bigl ( F(x) \bigr ) - F(x) \Bigr | \\
& \leq & \varepsilon + 2\varepsilon \cdot \max \|F_n\|_{\infty} +
3 \varepsilon + \varepsilon \leq 5  \varepsilon + 2\varepsilon
\cdot \max \|F_n\|_{\infty}
\end{eqnarray*}

\noindent Since the $F_n$'s are uniformly bounded, ${\widetilde
{F_n}}$ indeed converges uniformly to $F$ as $\varepsilon$ goes to
zero.
\end{proof}

\noindent As explained in Section $2$, the proof of
Claim~\ref{invariance} is based on a powerful result by
Katok~\cite{K} which is used for the proof of
lemma~\ref{approximation in measure lemma}.

\begin{proof}[{\bf Proof of Claim~\ref{invariance}:}]
Take $F \in L_{\infty}(M)$ 
and $\varphi$ a measure-preserving bijection on $M$. Consider a
sequence $\{F_n\}$ of uniformly bounded smooth functions which
converges in measure to $F$. Let $\varepsilon_n$ such that $F_n$
is an $\varepsilon_n$-approximation in measure of $F$. Choose
positive numbers $\delta_n$ so that $d(x,y) \leq \delta_n
\Rightarrow |F_n(x) - F_n(y)| \leq \varepsilon_n$. By repeatedly
using Lemma~\ref{approximation in measure lemma} we get
 a family of Hamiltonian diffeomorphisms
$\{g_n\}$
 such that
 $$ {\rm Vol}
\bigl ( \{ x \, | \, d \left (g_n(x),\varphi(x) \bigr)
> \delta_n \} \right )  \leq \varepsilon_n .$$
 Obviously
$$ \bigl |F_n \bigl (g_n(x) \bigr ) - F \bigl (\varphi(x) \bigr ) \bigr | \leq \bigl |F_n \bigl (g_n(x) \bigr ) -
F_n \bigl (\varphi(x) \bigr ) \bigr | + \bigl |F_n \bigl
(\varphi(x) \bigr) - F \bigl (\varphi(x) \bigr) \bigr| .$$ Our
choice of $\varepsilon_n$, $\delta_n$ and $g_n$ guarantees that
the above sum is smaller than $2\varepsilon_n$ outside a
$2\varepsilon_n$-measure exceptional set, and therefore that $\{
F_n \circ g_n \}$ converges in measure to $F \circ \varphi$. This
and the invariance of $\|\cdot\|'$ imply that
$$ ||| F \circ \varphi ||| \leq \liminf_n \|F_n \circ g_n
\|' = \liminf_n \|F_n \|' .$$

\noindent Since this is true for any sequence $\{F_n\}$ of
uniformly bounded smooth functions which converges in measure to
$F$, we conclude that $ ||| F \circ \varphi ||| \leq ||| F |||$.
Moreover, by applying the same argument to $F\circ\varphi$ and
$\varphi^{-1}$ we obtain that $ ||| F ||| \leq ||| F \circ \varphi
|||$, and the proof is complete.
\end{proof}

\section{Proof of Theorem~\ref{The_main_theorem}}

Let $\rho$ be an intrinsic bi-invariant pseudo-distance function
on Ham$(M,\omega)$ induced by some invariant norm on ${\cal A}$.
In order to determine whether $\rho$ is degenerate or not we will
use a criterion by Eliashberg and Polterovich~\cite{EP}. This
criterion is based on the following notion of ``displacement
energy" introduced by Hofer~\cite{H}.

\begin{definition}
For every open subset $A \subset M$ define its displacement energy
with respect to the pseudo-distance $\rho$ as
$$ e(A) = \inf \left \{ \rho({\mathbb \Id},\psi) \ | \ \psi \in {\rm Ham}(M,\omega),
\ \psi(A) \cap A = \emptyset \right \}, $$ and set $e(A) = \infty$
if the above set is empty.
\end{definition}

\begin{theorem}[Eliashberg-Polterovich] \label{Criterion_by_EP}
If $\rho$ is a genuine metric on Ham$(M,\omega)$ then the
displacement energy of every non-empty open set is strictly
positive.
\end{theorem}

\noindent This theorem allows to reduce the proof of
Theorem~\ref{The_main_theorem} to showing that the displacement
energy of some small ball vanishes. An argument borrowed
from~\cite{EP}, to be presented immediately below, further
reduces the problem to

\begin{claim} \label{small-support} If $\sup\{\|F_n\|_{\infty}\} < \infty$ and ${\rm
Vol} \bigl ({\rm Support}(F_n) \bigr) \rightarrow 0$, then $\| F_n
\| \rightarrow 0$.
\end{claim}

\noindent Indeed, choose an embedded open ball $B \subset M$ such
that its boundary $\partial B$ is an embedded sphere, and such
that there exists some Hamiltonian isotopy $\{ g_t \}$, $t \in
[0,1]$ generated by a Hamiltonian function $G(t,x)$ with $g_1(B)
\cap B = \emptyset$. Denote by $\Sigma_t $ the sphere
$g_t(\partial B)$. Consider the function $ K(t,x)$ obtained from
$G$ by smoothly cutting-off outside a neighborhood $U_t$ of
$\Sigma_t$. Note that the time-one-map of $K(t,x)$ also displace
$B$, i.e. $k_1(B) \cap B = \emptyset$. This is true since for
every $t \in [0,1]$ we have $k_t(\partial B) = g_t (\partial B)$.
Using Proposition~\ref{small-support}, we note that by decreasing
the sizes of the neighborhoods $U_t$ we can make the norm of
$K(t,x)$ arbitrary small. Hence the displacement energy of the
ball $B$ vanishes.

\noindent We are thus left with proving Claim~\ref{small-support}.
As explained in Section $2$, instead of proving it for
$\|\cdot\|$, we shall prove it for the extension $|||\cdot|||$
announced in Theorem~\ref{extension-to-L-infinity}.

\begin{proof}[{\bf Proof of Claim~\ref{small-support}}]

Let ${\mathbb \Id}_V$ stand for the characteristic function of the
set $V$. We first prove that
\begin{equation}\label{step-functions}
|||{\mathbb \Id}_V||| \rightarrow 0 \ {\rm as} \ {\rm Vol}(V)
\rightarrow 0 .
\end{equation}

\noindent Since $\| \cdot \|$ is not equivalent to $\| \cdot
\|_{\infty}$, and since $|||\cdot|||$ is an extension of
$\|\cdot\|$, for every $\varepsilon > 0$ there exists some
function $F \in {\cal A}$ with $\| F \| = |||F||| \leq
\varepsilon$, while $\| F \|_{\infty} = 1$. Assume that the
maximum of $F$ is obtained at some point $x_0 \in M$ and set $U$
to be a small-radius open set around $x_0$. Continuity allows us
to choose $U$ in such a way that $|F(x)|> 1 - \varepsilon$ for
every $x \in U$. Using the triangle inequality we obtain:
$$ |||  \langle F \rangle_U \cdot {\mathbb \Id}_U |||  \leq
 |||  \langle F \rangle_U \cdot {\mathbb \Id}_U  + \langle F
\rangle_{U^c} \cdot {\mathbb \Id}_{U^c} ||| + ||| \langle F
\rangle_{U^c} \cdot {\mathbb \Id}_{U^c}
 |||,$$
where $U^c = M \setminus U$. The left summand is estimated via
Lemma~\ref{claim-picewise-averaging}. To estimate the right
summand, recall that ${\rm Vol}(U)\langle F \rangle_U + {\rm
Vol}(U^c)\langle F \rangle_{U^c} = \langle F \rangle_M = 0$. We
therefore get:
$$ |||  \langle F \rangle_U \cdot {\mathbb \Id}_U |||  \leq  |||F||| + ||| {\frac { \langle F \rangle_U \cdot
 {\rm Vol}(U)} {{\rm Vol}(U^c)}} \cdot {\mathbb \Id}_{U^c} ||| .$$
Now, since $|||\cdot||| \leq C \|\cdot\|_{\infty}$, and since the
norm $\|{\frac { \langle F \rangle_U \cdot
 {\rm Vol}(U)} {{\rm Vol}(U^c)}} \cdot {\mathbb \Id}_{U^c}\|_{\infty}$ goes to zero with ${\rm Vol}(U)$,
 for $U$ with small enough measure we get
 $$  |||  \langle F \rangle_U \cdot {\mathbb \Id}_U |||
\leq |||F||| + \varepsilon \leq 2\varepsilon.$$ Due to the fact
that $ |\langle F \rangle_U| > 1-\varepsilon$, taking $\varepsilon
< 1/2$ we get $ ||| {\mathbb \Id}_U |||  < 4 \varepsilon$. Since
$|||\cdot|||$ is invariant under measure preserving bijections,
this applies to every set $V$ with the same measure as $U$ .

\noindent Now we can complete the proof of the Claim. Let $F \in
C^{\infty}(M)$ be supported in some compact set $U \subset M$ with
measure $\varepsilon$. Consider a finite partition of $U$ into
measurable sets $\{S_i\}_{i=1}^N$ with radius so small that
uniform continuity affirms $\max(F|_{S_i}) - \min(F|_{S_i}) \leq
\varepsilon$ for every $1 \leq i \leq N$. We have
$$|||F||| = ||| \sum_{i=1}^N F \cdot {\mathbb \Id}_{S_i}
||| \leq  ||| \sum_{i=1}^N \bigl  ( F - F(\eta_i) \bigr ) \cdot
{\mathbb \Id}_{S_i}
 ||| + |||  \sum_{i=1}^N  F(\eta_i) \cdot {\mathbb \Id}_{S_i}  |||,$$
 where $\eta_i$ is an arbitrary point in $S_i$. Without loss of
 generality we assume that $F(\eta_1) \leq F(\eta_2) \leq \ldots \leq F(\eta_N)$.
 Using the fact that $|||\cdot||| \leq C\|\cdot\|_{\infty}$ and the choice of $S_i$'s we get
$$ ||| F |||  \leq  C\| \sum_{i=1}^N \bigl  ( F - F(\eta_i) \bigr ) \cdot {\mathbb \Id}_{S_i}
 \|_{\infty} + |||  \sum_{i=1}^N  F(\eta_i) \cdot {\mathbb \Id}_{S_i}  ||| \leq
 C\varepsilon + |||  \sum_{i=1}^N  F(\eta_i) \cdot {\mathbb \Id}_{S_i}  |||$$
Next, in order to bound the last term on the right, we use Abel's
summation trick
$$  ||| \sum_{i=1}^N  F(\eta_i) \cdot {\mathbb \Id}_{S_i} ||| =
 ||| \sum_{i=1}^N  \bigl ( F(\eta_i) - F(\eta_{i-1}) \bigr ) \cdot {\mathbb \Id}_{\ \bigcup_{k=i}^{N} S_k} |||,$$
  where $F(\eta_0)$ is defined to be zero. Substituting this in the above inequality we conclude
$$ |||F|||  \leq  C\varepsilon +  \Bigl ( \sum_{i=1}^N  F(\eta_i) -
F(\eta_{i-1}) \Bigr ) \cdot  \max_i ||| {\mathbb \Id}_{\
\bigcup_{k=i}^{N} S_k} ||| \leq C\varepsilon + 2\|F\|_{_{\infty}}
\cdot \max_i |||  {\mathbb \Id}_{\ \bigcup_{k=i}^{N} S_k} |||.$$
Applying this estimate to a sequence of functions as in the
statement of the claim, recalling that $\varepsilon = {\rm
Vol}(\bigcup_{k=1}^{N} S_k)$ is the volume of the support, and
relying on (\ref{step-functions}), the proof of the claim is
complete.

\end{proof}

\section{Lemmas} \label{Section_Lemmas}

Here we prove Lemma~\ref{approximation in measure lemma},
Lemma~\ref{claim-picewise-averaging} and
Lemma~\ref{good_partition_lemma}. Recall that $M$ is a closed
connected symplectic manifold and $d$ is some Riemannian metric on
$M$.

\begin{proof}[{\bf Proof of Lemma~\ref{approximation in measure lemma}:}]
Fix $\varepsilon > 0$. Let $\{A_n\}_{i=1}^{N}$ be a family of
compact disjoint sets such that
\begin{enumerate}
    \item The diameter of each set $A_i$ is at most $\varepsilon$.
    \item ${\rm Vol}(\cup_{i=1}^N A_i) \geq {\rm Vol}(M) -
    \varepsilon$
\end{enumerate}
\noindent Next, let ${\widetilde B_i} = \varphi^{-1}(A_i)$, and
let $B_i$ be compact subsets of ${\widetilde B_i}$, such that
$\sum_{i=1}^N {\rm Vol}({\widetilde B_i} \setminus B_i) \leq
\varepsilon$. From Katok's Basic Lemma~\cite{K} we get a
Hamiltonian diffeomorphism $g$ satisfies $$\sum_{i=1}^N {\rm
Vol}(g(B_i) \triangle A_i) \leq \varepsilon.$$

\noindent We claim that $g$ is a good approximation in measure of
$\varphi$. To see that, let $C_i = \{x \in B_i \, | \, g(x) \in
A_i \}$ and denote by $C = \cup_{i=1}^N C_i$. Note that
$$ {\rm Vol}(C) = \sum_{i=1}^N {\rm Vol}(C_i) \geq \sum_{i=1}^N {\rm
Vol}(B_i) - \varepsilon \geq \sum_{i=1}^N {\rm Vol}(A_i) -
2\varepsilon \geq {\rm vol}(M) - 3\varepsilon .$$ Moreover, for
every $x\in C$ the points $g(x)$ and $\varphi(x)$ belong to the
same $A_i$. Since the diameter of the set $A_i$ is at most
$\varepsilon$, we conclude that $g$ is a
$3\varepsilon$-approximation in measure of $\varphi$.
\end{proof}

\begin{proof}[{\bf Proof of Lemma~\ref{claim-picewise-averaging}:}]
 In order to keep notation simple, let's assume we have only
two parts, $S_1$ and $S_2$. Fix $\varepsilon > 0$. Let
$\{U_j^1\}_{j=1}^{J_1}$ and $\{U_j^2\}_{j=1}^{J_2}$ be measurable
partitions of $S_1$ and $S_2$ respectively, such that:

\begin{enumerate}
\item All $U_j^1$'s have the same measure, and all $U_j^2$'s have
the same measure \item Inside all $U_j^i$'s the function $F$ does
not oscillate by more than $\varepsilon$
\end{enumerate}
Choose $\varphi_{j,k}^i$ to be measure preserving
transformations, not necessarily continuous, which map $U_j^i$
onto $U_k^i$. For every permutation $\pi$ of the set
$\{1,\ldots,J_i\}$, define $\varphi^i_{\pi}(x) =
\varphi^i_{j,\pi(j)}(x)$ if $x \in U_j^i$, and
$\varphi^i_{\pi}(x) = x$ if $x \not\in S_i$. Finally, define the
measure preserving bijections $\varphi_{\pi,\sigma} =
\varphi^1_{\pi}\circ\varphi^2_{\sigma}$. Since $|||\cdot|||$ is
invariant under measure preserving bijections, the triangle
inequality yields
$$|||\frac{1}{J_1! J_2!} \sum_{\pi,\sigma}
F\circ\varphi_{\pi,\sigma}||| \leq |||F|||.$$

\noindent Now, our choice of $U_j^i$ is such that for every $x \in
U_j^i$ we have $|\langle F \rangle_{U_j^i} - F(x)| \leq
\varepsilon$. Together with the equal measures of the $U_j^i$'s,
this means that for $x \in U_k^i$ we get $$|\frac{1}{J_i}
\sum_{j=1}^{J_i} F\circ\varphi^i_{k,j}(x) - \langle F
\rangle_{S_i}| \leq \varepsilon.$$ We thus infer the inequality
$$\| \frac{1}{J_1! J_2!} \sum_{\pi,\sigma}
(F\circ\varphi_{\pi,\sigma}) - (\langle F \rangle_{S_1} {\mathbb
\Id}_{S_1} + \langle F \rangle_{S_2} {\mathbb
\Id}_{S_2})\|_{\infty} \leq \varepsilon.$$ Since $|||\cdot||| \leq
C\|\cdot\|_{\infty}$, taking $\varepsilon$ to zero concludes the
proof.
\end{proof}

\begin{remark}\label{non-degenracy} {\rm
If $F$ is not continuous, then the choice of equal measure
$U^i_j$'s where $F$ has small oscillations may not be possible.
The argument, however, can be easily adapted to include the
non-continuous case as well. Lemma~\ref{claim-picewise-averaging}
is also the key behind the proof that $|||\cdot|||$ is indeed a
norm, namely that it vanishes only on the zero function. Indeed,
if $|||F|||$ were zero for a non-zero $F$, a piecewise averaging
of $F$ would generate a non-zero step function with vanishing
norm. Then, further piecewise averagings may be used to produce a
sequence of zero norm step functions, which converge uniformly to
a non-zero smooth function. This would mean that the original
norm, which coincides with $|||\cdot|||$ on smooth functions, was
already only a seminorm. We omit the details, because our main
results still hold even if $|||\cdot|||$ were only a seminorm. }
\end{remark}

\begin{proof}[{\bf Proof of Lemma~\ref{good_partition_lemma}:}]
According to Whitney's embedding theorem it follows that there
exists a smooth embedding $\Psi : M \rightarrow R^N$ for some
(large enough) $N$. Next, for $\alpha, \beta \in {\mathbb R}$ set
$\alpha {\mathbb K} + {\beta} = \{x \in {\mathbb R}^N \, | \,
\exists \ i \ {\rm such \ that} \  {x_i}-\beta \in \alpha {\mathbb
Z} \}$. Roughly speaking,  $\alpha {\mathbb K} + {\beta}$ stands
for the homothetic image of the ``standard" grid in ${\mathbb
R}^N$ translated in the direction of the vector $(1,\ldots,1)$.
For every $1 \leq j \leq J$, let $G_j$ be the
$\frac{\alpha}{4J}$-extension of the grid $\alpha {\mathbb K} +
{\frac{\alpha j}{J}}$. Set $\{V_i^j\}_{i=1}^{L_j}$ to be the
connected components of $\Psi(M) \cap (G_j)^c$. By choosing the
embedding coordinate-functions $\Psi_i$ to be Morse functions, we
can guarantee that the number of connected components is indeed
finite. It may well be the case that for a given $j$ some
$V_i^j$'s are zero-distance apart, but since our coordinates are
Morse functions, arbitrarily small translations of $G_j$ suffice
to guarantee positive distance separation between all $V_i^j$'s.

\noindent Now set $U_i^j = \Psi^{-1}(V_i^j)$. The first property
in the statement follows from the positive distance between the
$V_i^j$'s. Compactness guarantees that a small enough $\alpha$
implies the second property. The last property follows from the
fact that, regardless of $J$, the intersection of any $N+1$
different extended grids ($G_j$'s) is empty. Taking $J$ such that
$\frac{N+1}{J} < \varepsilon$ we are done.
\end{proof}

\section{Further Information Concerning our Norms}

\noindent Let $ \| \cdot \|$ be a Ham$(M,\omega)$-invariant norm
on ${\cal A}$ such that $\| \cdot \| \leq C\| \cdot \|_{\infty}$
for some constant $C$. Let $|||\cdot|||$ be the extensions of
$\|\cdot\|$ to $L_{\infty}(M)$  constructed in the proof of
Theorem~\ref{extension-to-L-infinity}. The main objective of this
section is to place the normed spaces $\bigl ( {\cal A},\|\cdot\|
\bigr )$ and $ \bigl (L_{\infty}(M),|||\cdot||| \bigr )$ in the
context of Banach (i.e. topologically complete) spaces of
functions, so that existing knowledge from this field be made
applicable to our context. For this purpose, we need to be able
to view our spaces as subspaces of Banach spaces of functions. It
can be easily seen that if the original norm $\|\cdot\|$ is
equivalent to $\|\cdot\|_{\infty}$, so is $|||\cdot|||$, and we
are in a Banach space setting. In the non-equivalent case we claim
the following. Here $\|\cdot\|'$ is the extensions of $\|\cdot\|$
to $C^{\infty}(M)$ constructed in the proof of
Theorem~\ref{extension-to-L-infinity}.

\begin{proposition}
Let $ \| \cdot \|$ be a Ham$(M,\omega)$-invariant norm on ${\cal
A}$ which is dominated from above by $\| \cdot \|_{\infty}$, but
not equivalent to it. The space $\bigl (L_{\infty}(M),|||\cdot|||
\bigr )$ then coincides with a dense subspace of the completion of
$\bigl (C^{\infty}(M),\|\cdot\|' \bigr )$. Moreover, this
completion can be viewed as a dense subspace of the space
$L_1(M)$ of integrable measurable functions on $M$, equipped with
a norm which is invariant under measure preserving bijections.
\end{proposition}

\begin{proof}[{\bf Sketch of the proof:}]
To establish the relation between $|||\cdot|||$ and the completion
of $||\cdot||'$, we need to show that if $\{F_n\}$ is a sequence
of uniformly bounded smooth functions tending in measure to $F$,
then $\{F_n\}$ is a Cauchy sequence in $\|\cdot\|'$ (which is
equivalent to showing that it is a Cauchy sequence in
$|||\cdot|||$). Indeed, let $F_n$ and $F_m$ both
$\varepsilon$-approximate $F$ in measure for some arbitrary small
$\varepsilon$. We can then write $F_n-F_m = G_{n,m}+H_{n,m}$,
where $G_{n,m}$ and $H_{n,m}$ are smooth and uniformly bounded,
$\|G_{n,m}\|_{\infty} \leq 2\varepsilon$, and the measure of the
support of $H_{n,m}$ is at most $2\varepsilon$.
Claim~\ref{small-support} now proves that $|||F_n-F_m|||
\rightarrow 0$ as $n,m \rightarrow \infty$.

\noindent We turn now to the second part of the proposition. First
we claim that there exist some constant $C$ such that $|||F|||
\geq C\|F\|_{L_1}$ for any essentially bounded measurable function
$F$. Indeed, set $M_F$ to be the median of $F$. Without loss of
generality we may assume that $M_F \geq 0$. Let $\{ x\in M \, | \,
F(x)
> M_F\} \subseteq A \subseteq \{x\in M \, | \, F(x) \geq M_F\}$, such
that ${\rm Vol}(A)=\frac{1}{2}$. Finally, let $B = \{ x\in M \, |
\, F(x) \geq 0\}$. We will argue under the assumption that $F$ is
zero-mean; the extension to general $F$ involves adding just one
dimension to our space of functions, and therefore follows
immediately. By Lemma~\ref{claim-picewise-averaging} we obtain
that $$|||F||| \geq |||\langle F \rangle_A {\mathbb \Id}_A +
\langle F \rangle_{A^c} {\mathbb \Id}_{A^c}||| = \langle F
\rangle_A|||{\mathbb \Id}_A - {\mathbb \Id}_{A^c}|||.$$ We
clearly have
$$ \langle F \rangle_A = 2 \int_A F(x) \, \omega^n \geq 2 \int_{B-A}
F(x) \, \omega^n ,$$ and therefore
$$ \langle F \rangle_A \geq \int_B F(x) \, \omega^n =
\frac{1}{2}||F||_{L_1}.$$
Together with the above estimate for
$|||F|||$, this yields $|||F||| \geq
\frac{1}{2}||F||_{L_1}\cdot|||{\mathbb \Id}_A - {\mathbb
\Id}_{A^c}|||$. The invariance property of the norm $||| \cdot
|||$ implies that the value of $|||{\mathbb \Id}_A - {\mathbb
\Id}_{A^c}|||$ depends only on the fixed ${\rm
Vol}(A)=\frac{1}{2}$. Thus, the inequality is proved.

\noindent Now, every Cauchy sequence of smooth functions in
$\|\cdot\|'$ (not necessarily uniformly bounded) is also a Cauchy
sequence in $L_1$, and is therefore convergent in measure. In
order to regard this limit in measure as an element of the
completion of $\|\cdot\|'$, we need to show that if two Cauchy
sequences in $\|\cdot\|'$, $\{F_n\}$ and $\{G_n\}$, converge in
measure to the same $F$, then both sequences have the same limit
in the completion, namely $\|F_n - G_n\|' \rightarrow 0$. For this
purpose set $H_n = F_n-G_n$. By definition, $H_n$ is a Cauchy
sequence in $\|\cdot\|'$ converging in $L_1$ and in measure to the
zero function. We need to prove that $\|H_n\|' \rightarrow 0$.
Again, we will carry out the proof in $|||\cdot|||$. By taking
small uniform perturbations we may also assume that ${\rm
Vol}(A_n) \rightarrow 0$, where $A_n = {\rm Support}(H_n)$. Next,
by applying a slight variant of
Lemma~\ref{claim-picewise-averaging} we get
$$|||H_n-H_m||| \geq |||(H_n-H_m){\mathbb \Id}_{A_m} + \langle
H_n-H_m \rangle_{A_m^c} {\mathbb \Id}_{A_m^c}|||.$$ Since $\langle
H_n-H_m \rangle_{A_m^c} \leq \|H_n-H_m\|_{L_1}\cdot{\rm
Vol}(A_m^c)$, this term goes to zero. We therefore conclude that
$|||(H_n-H_m){\mathbb \Id}_{A_m}|||=|||H_m - H_n\cdot{\mathbb
\Id}_{A_m}|||$ converges to zero as $n,m$ increase. Due to
Claim~\ref{small-support}, for every fixed $n$ the term
$|||H_n\cdot{\mathbb \Id}_{A_m}|||$ tends to zero with $m$. We
therefore conclude, as announced, that $|||H_m||| \rightarrow 0$.

\noindent Since we already know that $|||\cdot|||$ is invariant
under measure preserving bijections, the proof that the
completion is also invariant is straightforward. Note that since
we assume $\|\cdot\|$ is dominated by $\|\cdot\|_{\infty}$, but
not equivalent to it, Banach's Open Map Theorem implies that the
completion of $\|\cdot\|$ must in fact exceed the space of
essentially bounded measurable functions. The proof is now
complete.
\end{proof}

\noindent The literature contains much information concerning a
special subclass of the class of Banach norms on spaces of
functions, which are invariant under measure preserving
bijections. This is the subclass of the so called
\emph{Rearrangement Invariant Function Spaces}. The main (but not
only!) requirement is that the norm be monotone with respect to
the natural partial order on non-negative functions. Since an
explicit formulation will drag us into a long list of definitions
which are not relevant for this paper, we will make do here with a
reference. \cite{BS} introduces Rearrangement Invariant Function
Spaces in Chapter 2, Definition 1.4 (which relies on Definitions
1.1 and 1.3 from Chapter 1). The main classification results are
announced in Chapter 2, Theorem 5.15 and in Chapter 3, Theorem
2.12. Another thorough analysis from a somewhat different point of
view is available in \cite{LT} (Definitions 1.b.17 and 2.a.1, and
the results of the second section).

\noindent We cannot rule out the possibility that all normed
spaces $({\cal A},\|\cdot\|)$, which are invariant under
Hamiltonian diffeomorphisms, can be viewed as subspaces of
Rearrangement Invariant Function Spaces. The following example,
while not relating directly to the issue under discussion, serves
to indicate the kind of pathologies one might expect from norms
outside this class. Take the space ${\cal A} \bigoplus {\cal B}$,
where ${\cal B}$ is the space of functions on $M$ which attain
only finitely many values. It is straightforward to see that the
sum is indeed an algebraically direct sum. For an element $a+b \in
{\cal A} \bigoplus {\cal B}$ consider the functional
$\|a+b\|=\|a\|_1+\|b\|_{\infty}$. It is easy to check that
$\|\cdot\|$ is a norm invariant under measure preserving
diffeomorphisms, but not under measure preserving bijections (we
restrict our attention, of course, to measure preserving
bijections which keep the `rearranged' function inside our space).
It is also not hard to see that this norm is not bounded by
$\|\cdot\|_{\infty}$.

\noindent {\bf Acknowledgment:} This work is a part of the first
author's Ph.D. thesis, being carried out under the supervision of
Professor Leonid Polterovich at Tel-Aviv University. The first
author would like to thank Professor Polterovich for his guidance,
encouragement and continual support. Both authors wish to thank
Professor Leonid Polterovich for a thorough review of the paper
and many valuable suggestions, and Professor Vitali Milman for
useful comments.

\bigskip

\noindent

\begin{tabular}{@{} l @{\ \ \ \ \ \ \ \ \ \,} l }
Yaron Ostrover  & Roy Wagner \\
  School of Mathematical Sciences &  Computer Science Department \\
 Tel Aviv University  &   Academic College of Tel Aviv -- Yaffo \\
 Tel Aviv 69978, Israel &  4 Antokolsky St., Tel Aviv 64044, Israel  \\
yaronost@post.tau.ac.il  &
rwagner@mta.ac.il\\
\end{tabular}

\end{document}